
\input amstex
\documentstyle{amsppt}
\topmatter
\title Projective varieties with many degenerate
subvarieties \endtitle

\author
Emilia Mezzetti \endauthor
\address Dipartimento di Scienze
Matematiche - Universit\`a di Trieste - Piazzale Europa 1 -
34100 Trieste - Italy \endaddress
\email MEZZETTE\@UNIV.TRIESTE.IT (bitnet);
UTS882::MEZZETTE  or UTS340::MEZZETTE (decnet)\endemail

 \subjclass14J10, 14M07, 14N05\endsubjclass

 \thanks Work done with support of MURST and
CNR of Italy. \endthanks \endtopmatter
\document \head
Introduction \endhead Let $X\subset \Bbb{P}^{N}$ be a
projective non-degenerate variety of dimension $n$.   Then
the hyperplane sections of $X$ form a family of dimension
$N$ of varieties of dimension $n   -  1$, each one
generating a $\Bbb{P}^{N   -  1}$. The sections of $X$ with
the codimension two linear subspaces $L$ form a family,
parametrized by the Grassmannian $\Bbb G (N   -  2,N)$; if
$L$ varies in a suitable open subset of $\Bbb G (N   -
2,N)$ then each section has dimension $n   -  2$
 and generates  $L$. Some particular
codimension two linear space may possibly have improper
intersection with $X$, so that $X$ could contain some
subvarieties of dimension $n   -  1$ generating a $\Bbb{P}^{N
-  2}$. More generally, $X$ could have some \lq\lq degenerate''
subvarieties of codimension $h, h\ge 1$, i.e. subvarieties of
dimension $n   -  h$, contained in a linear space of dimension
at most $N   -  h   -  1 $.

The problem we shall study can be roughly formulated as
follows: give a classification of irreducible (possibly
smooth) non--degenerate varieties of dimension $n$ in
$\Bbb{P}^{N}$ containing a \lq\lq large--dimensional
family'' of degenerate subvarieties.

To clarify the meaning to be given to the expression \lq\lq
large--dimensional family'', let us consider some classical
examples.
 \example{Example 1}
Let $n=2, N=4.$ It is well known (see for instance \cite
{B}) that the only surfaces of $\Bbb{P}^{4}$ containing
a $2$ -- dimensional family of conics are the projections of
 the
Veronese surface of $\Bbb{P}^{5}$. There are no irreducible
 non-degenerate surfaces in $\Bbb{P}^{4}$ with a $3$ --
dimensional family of conics. This result was improved by
Corrado Segre (\cite{S1}, see also \cite{CS}), who proved
that, if an irreducible surface of  $\Bbb{P}^{4}$ contains a
family of dimension $2$ of plane curves of degree $d$, then $d
\leq 2$. So either the curves are conics, and we fall in the
above case, or they are lines, and the surface is a plane.
\endexample \example{Example 2} The study of varieties
containing \lq\lq many'' linear subspaces is also classical: it
corresponds to the case $N=n+1$ in the problem above. In the
notes \cite{bS}, Beniamino Segre
stated the problem of: \roster
\item determining what integers $c > n  - k$ can be dimension
of a family of linear subspaces of dimension $k$ contained in
a variety of dimension $n$;
\item classifying varieties of dimension $n$ containing
such a family of subspaces.
\endroster
He proves that the maximal dimension is $\delta (n,k) =
(k+1)(n  - k)$, which is achieved only by linear varieties.
Moreover some classification results are given, but there is
no complete solution to the problem (see $\S 2$).
\endexample
\example{Example 3}
Let $n=2, N=5.$ In \cite{S1}  and \cite{S2} Corrado
Segre studied the surfaces $S$ in $\Bbb{P}^{5}$ with a family
$\Cal F$ of dimension $c \ge 2$ of curves of $\Bbb{P}^{3}$. His
result may be summarized as follows: \roster
\item if $c=3$, then $S$ is a rational normal scroll of degree
$4$;
\item if $c=2$, then there are two cases: either $S$ is any
surface contained in a $3$ -- dimensional rational normal
scroll of degree $3$ or in a cone over a Veronese surface, or
the curves in $\Cal F$ have degree at most $5$ and a finite
number of cases is possible.
\endroster \endexample
Now we give a more precise formulation of the general problem.
\subsubhead Problem \endsubsubhead
Let $2 \leq n < N, 1\leq h < n$ be integers. Classify
irreducible non--degenerate varieties of dimension $n$ in
$\Bbb{P}^{N}$, containing an algebraic family $\Cal F$ of
dimension $c \geq h+1$ of subvarieties of dimension $n - h$,
each one spanning a linear space of dimension at most $N - h -
1$.

In this paper, we propose an approach to this problem, which
extends the classical one and stems from the notion of $foci$
of
a family of varieties: a classical notion that was recently
rediscovered and used in various situations (geometry of
canonical curves
\cite {CS}, lifting
problems \cite{CC}, \cite{CCD}, \cite{Me}). We
assume that $c=h+1$ and that a general variety in $\Cal F$
generates a linear space of dimension exactly $N- h - 1$, and
study the foci of the $1^{st}, 2^{nd}$,... order of the family
of  $\Bbb{P}^{N - h - 1}$'s spanned by the varieties of $\Cal
F$. We get ($\S 1$) that the solutions may be  essentially
divided in two parts: varieties contained in a higher
dimensional variety with similar properties but bigger $h$, and
varieties that we call {\it of isolated type}, for which there
is a bound on the degree of the degenerate subvarieties. So one
sees that to completely solve the problem for varieties of
codimension $N - n$, one has to solve before the problem for
varieties of codimension $<N - n$.

Then we discuss the first particular cases of the question. In
$\S 2$ we collect the results of B. Segre for $N - n=1$ and
some recent related results by Lanteri--Palleschi. In $\S 3$ we
consider the case of varieties of codimension two in
$\Bbb{P}^{N}$: we treat the cases $h=n-1, n-2$, corresponding
to varieties containing many plane curves or surfaces of $\Bbb
P^3$. For $n=3$,
we see that the problem is easy if $h=1$, since it reduces, by
cutting with a general hyperplane, to the case $n=2$ which is
well known; if $h=2$, the problem for the non--isolated type
seems not to be easy. We discuss some examples but we have not
uniqueness results. In $\S 4$, we consider varieties with
$N-n=3$; since the very long proof of the theorem of
C. Segre, on surfaces in $\Bbb{P}^{5}$ containing a
$2$-dimensional family of space curves (see Example 3), seems
to us not competely correct, we have tried to rewrite and
translate it in a modern language, clarifying some rather
oscure points. Then we make some remarks on the next cases.

\subsubhead Notations \endsubsubhead
We work over $k$, an algebraically closed field of
characteristic $0$. $\Bbb P ^N$ will denote the projective
space of dimension $N$ over $k$.
By {\it variety} we will always mean an algebraic reduced
scheme over $k$. If $V$ is a subscheme of $W$,  $\Cal
N_{V,W}$  will denote the normal sheaf of $V$ in $W$, $T_W$
the tangent sheaf of $W$.
 \head 1. Foci of families of linear spaces
\endhead
In this section we will  first recall some  facts about
foci of families of projective varieties (see
\cite{S}, \cite{CS}, \cite{CC}).

Let us fix integers $N>0$ and  $h<N - 1$. Let $Z$ be a
non--singular quasi--projective variety of dimension $c$
and  $\Phi \subset Z\times \Bbb P ^N$ be a flat family of
projective irreducible subvarieties of codimension $h+1$ in
$\Bbb P ^N$. Let $q_1, q_2$ be the natural
projections from  $ Z\times \Bbb P ^N$ to $Z, \Bbb
P ^N$; $p_1, p_2$ their restrictions to $\Phi$. We will assume
that the natural map from $Z$ to the
Hilbert scheme of subvarieties of $\Bbb P ^N$ is finite. The
family $\Phi$ is said to be $non-degenerate$ if the image of
the
projection $p_2$ has dimension $N - h - 1 +c$, i.e. if $p_2$ is
generically finite onto its image. A {\it fundamental point}
for the family $\Phi$ is by definition a point $p$ of $\Bbb P
^N$ such that
$dim \ p_1(p_2^{- 1}(p))>N-h-1+c-dim \ p_2(\Phi )$; in
particular, if $\Phi$ is non--degenerate, then $p$ is a
fundamental point if it belongs to infinitely many varieties of
$\Phi$.

Let us consider now $\tilde \Phi$, a desingularization of
$\Phi$, that we may assume to be flat over $Z$, with smooth
irreducible fibers. Moreover we may assume that the
fibers of $\tilde \Phi$ are desingularizations of the
fibers of $\Phi$. Denote by $u$ the natural map $\tilde \Phi
@>>>  \Phi \subset Z\times \Bbb P ^N$: its differential $du:
T_{\tilde \Phi}\rightarrow u^*T_{  Z\times \Bbb P ^N}$ is an
injective morphism of sheaves, whose cokernel is $\Cal N$, the
normal sheaf to the map $u$, a non--necessarily torsion--free
sheaf. The following exact sequence of sheaves on $\tilde\Phi$
defines a locally free sheaf $T(p_2)$
 of rank $c$:$$0\rightarrow
T(p_2)\rightarrow u^*T _{ Z\times \Bbb P ^N}\rightarrow
u^*q_2^*T_{\Bbb P^N} \rightarrow 0.$$ We have the following
exact commutative diagram: $$\CD {} @.   0   @.  0
@.  {}  @. {} \\
 @.      @VVV   @VVV  @.  @.  \\
0  @>>>  ker \ \chi  @>>>  T(p_2)  @>\chi>>
 \Cal N  @.  {}  \\
@.      @VVV   @VVV  @VVV    @.  \\
0  @>>>  T_{\tilde \Phi}   @>>>  u^*T_{Z\times \Bbb P^N}
 @>>> \Cal N  @>>>  0  \\
@.   @V d(q_2u) VV   @VVV   @.   @.   \\
{}   @.  u^* q_2^*T_{\Bbb P ^N}  @=  u^*q_2^*T_{\Bbb P
^N}  @.  {}   @.     \\  @.  @. @VVV @.  @.  \\ {} @. {}
@. 0 @. {} @. {}\endCD$$ where $\chi$ is the global
characteristic map for $\tilde \Phi$; it is such that $$rk \
\chi=dim \ p_2(\Phi )-(N-h-1), \ rk(ker \ \chi)= dim \ \Phi
-dimp_2(\Phi ).$$ If $\Phi$ is non--degenerate, then $rk \
\chi=c$ and $rk(ker \ \chi)=0$.

Note the following facts (\cite{CC}):

 - if $\tilde H$ is a general fiber of $\tilde \Phi$, then
$\Cal N\mid _{\tilde H}$ is isomorphic to the
normal sheaf to the map $\tilde H\rightarrow\Bbb P^N$
induced by $\tilde \Phi$;

 - the restriction of
$\chi$ to a fiber $\tilde H=(uq_1)^{-1}(z), z\in Z$, is a
morphism
$\chi_{\tilde H}: T_{z,Z}\otimes\Cal O_{\tilde H}\simeq\Cal
O_{\tilde H}^c\rightarrow \Cal N_z$;

- on an open subset of $\tilde H$ $ker(\chi_{\tilde H})$
coincides with $(ker \ \chi)_{\tilde H}$.

A point $P$ of $\tilde \Phi$ (or its image $p$ in $\Bbb P^N$)
is called a {\it $1^{st}$ order focus} of $\Phi$ on $H$ if the
map $\chi _P: T(p_2)\otimes k(P)\rightarrow \Cal N\otimes
k(P)$ has rank less than $c$; it is called a \it cuspidal
point \ \rm if the map $T_{\tilde \Phi ,p}\rightarrow
u^*T_{Z\times\Bbb P^N,p}$ is not injective. The cuspidal
(focal) locus is the set of all cuspidal (focal) points.

Note that if $P$ is cuspidal, then $\Cal N$ is not free at
$P$ and $u(P)$ is singular both inside  $\Phi$ and in its
fiber of $p_2$. The focal locus in $\tilde H$ is closed off
the cuspidal locus; it is in fact the degeneration locus of
the map $\chi _{\tilde H}$.

       \proclaim{1.1. Lemma} Let
$\Phi$ be a family as above. Then any fundamental point of
$\Phi$ is a focus or a cuspidal point. \endproclaim
\demo{Proof} See \cite{CC}, Prop.1.7.\enddemo

 Assume now that $\Phi$ is a non--degenerate family of smooth
subvarieties $H$ of $\Bbb P ^N$ with $c=h+1$.  For any $H$  in
an open subset $Z_2$ of $Z$, we have the codimension one closed
subscheme $F_H^1$ of the $1^{st}$ order foci on $H$; so we may
consider the family $\Phi _2 \subset Z_2 \times \Bbb P ^N$ of
such subschemes. Assume that $F^1_H$ is integral if $H\in Z_2$:
by definition, the $1^{st}$ order foci of $\Phi _2$ on $F_H^1$
are called the {\it $2^{nd}$ order foci} of $\Phi$ on $H$; they
form a closed subscheme $F_H^2$ of $F_H^1$ off the cuspidal
points.

  And so on by
induction; fix an integer $i\ge 2$: then the scheme
$F_H^{i}$ of $i^{th}$ \it order foci \rm of $\Phi$ is by
definition the scheme of $1^{st}$ order foci of the family
$\Phi _i$ of schemes of $(i - 1)^{th}$ order foci. It is well
defined under the assumptions $F_H^{i-1}$ integral for $H$
general in an open subset $Z_i$ of $Z$ and $codim _H
F_H^{i-1}=i-1$ or, equivalently, $F_H^k\ne F_H^{k-1}$ for
$k\leq i - 1$.  Note that, by 1.1, if $p$ is a fundamental
point of $\Phi$, then either  it is a singular point of some
focal locus or it is a focus of any order for the family.

 Now we come to our problem.
Let $X \subset \Bbb P ^N$ be an integral non--degenerate
projective variety of dimension $n\geq 2$. Let $h$ be an
integer, $1\leq h <n$. Assume that there is a family $\Cal F
 \subset
Z\times X$  of subvarieties of dimension $n - h$ of X,
parametrized by an integral smooth variety $Z$ of dimension
$c= h+1$. We will say that $X$ contains a family of
dimension $c$ of (or  $\infty ^c$) subvarieties of
dimension $n - h$. Let us call $Y_z$, or simply $Y$, a general
element of the family. Assume that the linear span of $Y$,
$\langle Y \rangle$, has dimension  $N - h - 1$. We would like
to classify such kind of varieties.

We define now a family of linear spaces associated to
$\Cal F$: let $\Phi _1$ be the family of $\Bbb P^{N - h - 1}$'s
spanned by the varieties of the family $\Cal F$; it is
parametrized by the same variety $Z$ as $\Cal F$, so  it has
dimension $c$.

\proclaim{1.2  . Lemma} Let $X \subset \Bbb P ^N$
be an integral variety of dimension $n$
containing a family $\Cal F$ of dimension at least
$h+1$ of subvarieties of dimension $n - h$. Then the
points of X are fundamental points for the family $\Phi _1$ of
the linear spaces generated by the varieties of $\Cal F$.
\endproclaim \demo{Proof} Let $p$ be a point of $X$; passing
through $p$ imposes $h$ conditions to the varieties in $\Cal
F$, so the subfamily $\Cal F_p$ of varieties of $\Cal F$
through $p$ has dimension at least $1$. Hence $p$ is a
fundamental point for $\Cal F$ and also for $\Phi _1$ and the
lemma follows. \enddemo

\proclaim{1.3. Theorem} Let $X \subset \Bbb P ^N$
be an integral variety of dimension $n$
containing a family $\Cal F$ of dimension
$c=h+1$ of integral subvarieties of dimension $n - h$. Let $Y$
be a general variety of $\Cal F$ and assume that $Y$ spans a
$\Bbb P ^{N - h - 1}$. Then one of the following happens:
\roster
\item there exists an integer $r$, $1\leq r < N - n$, such
that $X$ is contained in a variety $V _r$ of dimension at
most $N - r$ containing $\infty ^{h+1}$ varieties of dimension
$N - h - r$, each one contained in a linear space of dimension
$N - h - 1$;
\item deg $Y$ is bounded by a function of $h, N-n$.
\endroster
\endproclaim
\demo{Proof}   Let $\Phi _1$ be the family of
linear spaces generated by the varieties of $\Cal F$
and $H=\langle Y \rangle$ be a general element of $\Phi _1$.
Then $F_H^1$, the locus of the $1^{st}$ order foci of $\Phi _1$
on $H$, is defined by the vanishing of  $det \ \chi _H$, where
$\chi _H$ is as follows: $$\CD T_{Z,z} \otimes
\Cal O_{H}   @>\chi _H>>  \Cal N_{H,\Bbb P ^{N}}
 \\
 @\vert            @\vert \\
\Cal O_{H}^{h+1}   @>\chi _H>>   \ \Cal O_{H}(1)^{h+1}. \endCD
$$
If $det \ \chi _H\equiv O$, then $\Phi _1$ is degenerate, so
$V_1:=\bigcup_{H\in \Phi _1}H$ is a variety of dimension $\le
N - 1$ containing $X$ which is covered by $\infty ^{h+1}$
linear spaces of dimension $N - h - 1$, and $(1)$ holds with
$r=1$. If $det \ \chi _H\not \equiv O$, the condition $det
\ \chi _H= O$ defines a hypersurface $F_H^1$
 of degree $h+1$ in
$H\simeq \Bbb P ^{N - h - 1}$ containing $Y$. If $N - n =1$,
only the first case may happen; if $N - n=2$ and the second
case happens,
 we have $dim \ Y=dim \ F_H^1$ so $deg \ Y\leq h+1$
and the theorem is proved.

If $N - n>2$, we  consider the family $\lbrace
F^1_H\rbrace$. If the general $F_H^1$ is integral, we set $\Phi
_2=\lbrace F^1_H\rbrace$; otherwise we replace
$F_H^1$ with one of its irreducible components
containing $Y$ with reduced structure. If $\Phi _2$
is degenerate, case (1) happens; if it is
non--degenerate,by (1.2) any point $p$ of $X$ is a fundamental
point of $\Phi _1$; for any $H$ containing it, $p$ is a focus
of $1^{st}$ order, so it is a fundamental point for $\Phi _2$:
by (1.1), either $p\in Sing \ F^1_H$, the singular locus of
$F^1_H$, or $p\in F^2_H$, the focal locus of $\Phi _2$ on $H$.
Therefore $Y\subset Sing \ F^1_H\cup F^2_H$; $Y$ being
integral, it is completely contained in one of them.

 If $N - n=3$
then $deg \ Y \leq max \lbrace deg \ Sing \ F^1_H, deg \
F^2_H\rbrace$; otherwise we proceed in this way, by defining
families $\Phi _i=\lbrace F_H^{i-1}\rbrace$ for $i\leq
N-n-1$. If one of these families is degenerate, (1) happens.
Otherwise $$deg \ Y\leq max\lbrace deg \ Sing \ F_H^{N-n-2},
deg \ F_H^{N-n-1}\rbrace$$.

To have an estimate of such bound, note that  $dim \
F_H^i=N-h-1-i$; moreover, denoting $d_i=deg F^i_H
(i=1,...,N-n-1), s_i=deg \ Sing F_H^i$, there are the
relations (\cite{CC}):

- $d_i\leq(N+1)d_{i-1}+deg \ K_{i-1}$ where $K_{i-1}$ is the
canonical class on $\tilde F_H^{i-1}$ (the desingularization
of $F_H^{i-1}$);

- $deg \ K_{i-1}\leq
(1-N+h+i)d_{i-1}-2s_{i-1}+k_i(k_i-1)(i-1)+2e_ik_i-2$;

- $s_{i-1}\leq \frac {k_i(k_i-1)} {2} (i-1)+e_ik_i
$

where $k_i=\left[\frac {d_{i-1}-1} {i-1}\right],
0\leq e_i=d_{i-1}-1-k_i(i-1)<i-1$.

 So we get the recursive formula:$$d_i \le \frac
{d_{i - 1}^2} {i - 1} +(h+2+i)d_{i - 1}+ o(d_{i - 1})$$
where the latter term is independent of $h$; we conclude that
$deg \ Y$ is bounded by a function of $h$ and $N-n$ having
the same order as $d_{N-n-1}$.
\enddemo
\definition{1.4.
Definition} If a variety $X$ satisfies condition $(2)$ of $1.3$
we say that $X$ is {\it of isolated type}. \enddefinition

\remark{1.5. Remarks}
\roster
\item"(i)" The hypersurfaces of degree $h+1$ of $\Bbb P^{N - h
- 1}$ arising as $1^{st}$ order foci are not general; they
are in fact linear determinantal varieties. For
example in $\Bbb P^3$ the general surface of
degree $d\geq 4$ is not determinantal (see \cite{H}). If $h=1$
and $N\geq 5$ they are quadrics of rank at most $4$.
 \item"(ii)" The codimension $2$ subvarieties of $\Bbb P^{N - h
- 1}$ which are $2^{nd}$ order foci of $\Phi _1$ are
aritmetically Cohen--Macaulay. In fact they are defined by
the condition $rk \ \chi _1\leq h$ where $\chi _1: \Cal
O^{h+1}_{F^1_H}\rightarrow \Cal O_{F^1_H}(1)^{h+1}\oplus
\Cal O_{F^1_H}(h+1)$. By considering a lifting of $\chi _1$
to $H$, we see that it degenerates along a subvariety of
codimension $2$ in $H$, of degree
$$deg \ c_2(\Cal
O_{H}(1)^{h+1}\oplus  \Cal O_{H}(h+1))= (h+1)(3h+2)/2.$$ This
variety is linked in a complete intersection of type $(h+1,
2h+1)$ to a variety of degree $h(h+1)/2$ defined by a $h \times
(h+1)$--matrix of linear forms. In particular, if $h=1$ it is a
Castelnuovo variety. \item"(iii)" When considering the
$i^{th}$ order foci of $\Phi _1$, $i\geq 3$, we find the normal
bundle $\Cal N_{F_H^{i - 1},H}$ which is not decomposable in
general, so we cannot lift the map $\chi _{i - 1}$ to a map
between bundles on $H$. For example, if $N=6, n=2, h=1$, the
$2^{nd}$ order foci  form a Castelnuovo surface $S$ of degree
$5$ in $\Bbb P^4$, whose normal bundle is not decomposable. In
this case a bound on the degree of $Y$ is given by deg $
c_1(\Cal N_{S\mid \Bbb P^4}) = 5 (deg S)^2 +2\pi - 2 - deg S =
22$ (here $\pi$ denotes the sectional genus of $S$).
\item"(iv)" If $h>1$, by the proof of (1.3) we have that $d_1
\le h+1$, and $d_i$ is bounded by a function of
order $\frac {d_{i - 1}^2} {i - 1}$. So we get an upper bound
of order $$\frac {h^{2^{i - 1}}} {(i - 1)(i - 2)^2 (i -
3)^{2^2}...2^{2^{i - 3}}}.$$   \item"(v)" Note that if $V_r$ is
a variety as in (1), then any variety $X$ of dimension $n$
contained in $V_r$ satisfies the assumption of the Theorem.
\endroster \endremark

In the next sections we will treat the first particular cases
of the question, for $N - n=1,2,3$. As we will see, for
fixed codimension, the situation becomes simpler as $h$
decreases.

\head 2. The case $N - n=1$ \endhead
In this section we collect the known results about the
codimension $1$ case, i.e. the case of a variety $X$ of
dimension $n$ embedded in $\Bbb P^{n+1}$ and
 containing a family of dimension $c\geq {h+1}$ of $\Bbb
P^{k}$'s, $k=n - h$.

In the above quoted papers  \cite{bS} B. Segre, who does not
make restrictions on the dimension of the projective space
containing $X$,
 first states some existence results. Assume $\delta = \delta
(n,k) = (k+1)(n - k)$, the dimension of the Grassmannian of
$k$--planes in $\Bbb P^n$ $\Bbb G (k,n)$. He proves that:
\roster
\item"--" if $k=1$, then $c$ may assume all values from $n$ to
$\delta (n,1)$;
\item"--" if $k=2$, then $c$ may assume all values from $n-1$ to
$\delta (n,2)$, except $\delta -1$; we say that $c= \delta -1$
is a \it gap\rm;  \item"--" if $k>2$, then the following $k - 1$
numbers are gaps for $c$: $\delta - k+1,...,\delta - 1$;
$\delta - k$ is not a gap; then there are $k - 2$ gaps: $\delta
- 2k+2,...,\delta - k - 1.$ \endroster
As for uniqueness
results, he proves: \roster
\item if $c=\delta$, then $X$ is linear;
\item if $c=\delta - k$, then $X$ is a scroll in $\Bbb P^{n -
1}$'s, or a quadric if $k=1$;
\item if $c=\delta - 2k+1$, then $X$ is a quadric.
\endroster

There is one more case he treats, i.e. $k=1, c=\delta - 2$. He
 quotes results by Togliatti and Bompiani (\cite{T}, \cite{Bo})
claiming that in this case $X$ has to be a scroll in $\Bbb
P^{n - 2}$'s, or a quadric bundle ($n\geq 4$), or a section of
$\Bbb G(1,4)$ with a $\Bbb P^7$ (the latter variety is not a
hypersurface). Recently, this claim
has been proved by Lanteri--Palleschi (\cite{LP}) under
smoothness assumption, by adjunction--theoretic techniques.
Unfortunately, it is not clear if this result is true without
assuming that $X$ is smooth. There is work in progress on this
subject by E.Rogora.

\head 3. The case $N - n=2$ \endhead
In this section $X$ is a codimension $2$ subvariety of $\Bbb
P^N$ containing $\infty ^{h+1}$ varieties of dimension $n - h$
each one generating a $\Bbb P^{N - h - 1}$. Then, by
$(1.3)$, there are the following possibilities:
\roster
\item $X$ is contained in a hypersurface $V_1$ containing
all the linear spaces of the family $\Phi _1$;
\item $deg \ Y \leq h+1$. \endroster
Let us study the first particular cases.\smallpagebreak

$\boxed{a) \  h=n-1}$

The varieties $Y$ are plane
curves. There are two subcases:
\subsubhead\nofrills{a1) } \endsubsubhead $X$ is contained in a
variety $V_1$ of dimension $n+1$ containing at
least $\infty ^n$
planes. By $\S 2$ these varieties are classified only for
$n\leq 3$: if $n=2$ they are linear, which is impossible in our
case because $X$ is non--degenerate; if $n=3$, then $V_1$ is a
quadric or a scroll in $\Bbb P^3$'s;
  \subsubhead\nofrills{a2) } \endsubsubhead the varieties $Y$
are plane curves of degree $d\leq n$, $d>1$.

If $n=2$, we find again the classical case of surfaces of $\Bbb
P^4$ containing $\infty ^2$ conics (see the Introduction).

If $n=3$, the curves are conics or cubics. Looking at the
list of the known smooth threefolds of $\Bbb P^5$, we find two
examples, precisely the Bordiga and
Palatini scrolls (see \cite{O}). The Bordiga scroll is an
arithmetically Cohen--Macaulay $3$--fold in $\Bbb P^5$ which is
defined by the $(3\times 3)$--minors of a general $(4 \times
3)$--matrix with linear entries. Its general hyperplane section
is a Bordiga surface of $\Bbb P^4$, a rational surface $S$ of
degree $6$, which can be realized as the blowing--up of $\Bbb
P^2$ in $10$ points, embedded by the complete linear series of
quartics through these points. Then any line through $2$ of the
$10$ points corresponds to a conic on $S$ and any cubic through
$9$ of them corresponds to an elliptic cubic; so $S$ contains
$45$ conics and $10$ plane cubics. Hence the Bordiga scroll
contains $\infty ^3$ conics and $\infty ^3$ plane cubics.

The Palatini scroll $X$ is an arithmetically
Buchsbaum $3$--fold
whose ideal has the following $\Omega$--resolution (see
\cite{Ch}):$$
0 \to 4\Cal O_{\Bbb{P}^{5}}  \to \Omega_{\Bbb{P}^{5}}^1(2)
\to \Cal I_X(4) \to 0;\tag 3.1$$
i.e. $X$ is the degeneration locus of a bundle map defined by
$4$ independent sections of $\Omega_{\Bbb{P}^{5}}^1(2)$. A
general hyperplane section of $X$ is a rational surface $S$ of
degree $7$ and sectional genus $4$ (see \cite{Ok});
the map that gives the embedding in $\Bbb P^4$ is associated
to the linear system of the sextics plane curves with $6$ double
points and $5$ simple points in assigned general position. Any
line through two of the $6$ double points corresponds to a
conic on $S$ and any cubic through the $6$ double points and
$3$ of the simple points corresponds to an elliptic cubic. So
also the Palatini scroll contains $\infty ^3$ conics and
$\infty ^3$ plane cubics.

Let us give another construction of such plane curves on $X$
by means of a geometrical interpretation of the exact sequence
(3.1); it is the direct generalization of a construction by
Castelnuovo which gives the Veronese variety in $\Bbb P^4$(see
\cite{Ca}, \cite{O}).
 Let us define $V=H^0(\Cal O_{\Bbb P^5}(1))^*$;
by the Euler's sequence, $H^0(\Omega ^1(2)) \simeq \Lambda
^2V^*$, so any section of $\Omega^1(2)$ may be seen as a
bilinear antisymmetric form on $V^*$, associated to an
antisymmetric $(6\times 6)$--matrix $A=(a_{ij})$,
$i,j=0,...,5$, with constant entries. It defines a null
correlation $\Phi: \Bbb P^5 \to \check {\Bbb P^5}$, a rational
linear map which associates to a point $p$ a hyperplane through
$p$. $\Phi$ associates to a line $l$ a $\Bbb P^3$, the
intersection of the hyperplanes corresponding to the points of
$l$. Let us consider the set of lines $l$ such that $l \subset
\Phi (l)$: it is classical that this is a linear complex
$\Gamma$ of lines in $\Bbb P^5$, i.e. a hyperplane  section of
the Grassmannian $\Bbb G(1,5)$, which has equation precisely
$\sum_{i,j}a_{ij}p_{ij}=0$ ($p_{ij}$ coordinates of a line in
the Pl\"ucker embedding of  $\Bbb
G(1,5)$).

Four sections of $\Omega ^1(2)$
define  four matrices $A_1,...,A_4$, $A_k=(a_{ij}^k)$, four null
correlations $\Phi _1,...,\Phi_4$, four linear complexes
$\Gamma _1,... \Gamma _4$; for a general point $p$, the four
corresponding hyperplanes via $\Phi _1, ... ,\Phi _4$ intersect
along a line: by definition $X$ is the set of points $p$ such
that they intersect along a plane $\pi$, i.e. the matrix $$M:=
\left( \matrix  \sum a_{0i}^1x_i  &   \hdots & \sum a_{5i}^1x_i
 \\
\hdots  & \hdots  &  \hdots  \\
\sum a_{0i}^4x_i  &   \hdots &
\sum a_{5i}^4x_i \endmatrix \right) $$ has rank $<4$. So $X$ is
defined by a system of equations of degree $4$.

Let us recall that a linear complex $\Gamma$ of lines of $\Bbb
P^5$ is said to be non-special if the associated matrix $A$ is
non--degenerate; special if it is degenerate. The rank of $A$,
$\rho (A)$, is always an even number; if $\rho (A)=4$, there
is a line, the center of $\Gamma$, which is met by all lines
in $\Gamma$; if $\rho (A)=2$, there is a $\Bbb P^3$ which is
the center of the complex. Given four linear complexes $\Gamma
_1,..., \Gamma _4$, with associated matrices $A_1...,A_4$,
they define a linear system of complexes; the special
complexes of the system correspond to the zeros $(\lambda
_1,...,\lambda_4)$ of the pfaffian of the matrix $\lambda
_1A_1+...+\lambda _4 A_4$, so they are parametrized by a cubic
surface of $\Bbb P^3$. The union of the lines, centers of these
complexes, is the degeneracy  locus of the bundle map
 $$4\Cal O_{\Bbb{P}^{5}}
\to \Omega_{\Bbb{P}^{5}}^1(2)
$$ given by $A_1...,A_4$, i.e. the Palatini scroll. The
intersection of $\Gamma
_1,..., \Gamma _4$ is a family of lines of dimension $4$,
$\Sigma _4$, such that there is one line of $\Sigma _4$
through a general point $p$ of $\Bbb P^5$; but if $p \in X$,
there is a pencil of lines of $\Sigma _4$ through $p$,
spanning a plane $\pi$. One can easily see, by exactly the
same argument as in \cite{Ca},
that $X \cap \pi$ is the union of
the point $p$ and a plane cubic, generally not containing $p$.
So we get a family of plane cubics on $X$ parametrized by the
points of $X$.

Assume now that one of the matrices  $A_i$, $1\leq
i\leq 4$, has rank $2$; in particular,
we may suppose that $A_1$
has the following canonical form: $$A_1=\left( \matrix
0 & 1 & 0 & \hdots & 0\\
-1 & 0 & 0 & \hdots & 0\\
0 & 0 & 0 & \hdots & 0\\
\vdots & \vdots & \vdots & \vdots & \vdots
\endmatrix \right);$$
in this case the center of the complex $\Gamma _1$ is the
space $H: x_0=x_1=0$. Note that $M$ takes the form $$M=
\left( \matrix \matrix x_1 & -x_0 \endmatrix & & \matrix 0 & 0
& 0 & 0 \endmatrix  \\
 & &\\
        \ P & & N\endmatrix
\right).$$ It easily follows that $X$
  splits; in fact its components are $H$ and the variety $Y$
defined by the maximal minors of the $3\times 4$ matrix  $N$;
hence $Y$ is a Bordiga scroll. So this construction gives both
the known examples of smooth threefolds with $\infty ^3$ plane
curves, of isolated type.

The above construction may be generalized to higher
dimension, giving examples of (singular) varieties of isolated
type with $N - n=2, h=n - 1$, for any $n\geq 3$.\smallpagebreak

$\boxed{b) \  h=n-2}$

The varieties $Y$ are
surfaces of $\Bbb P^3$.
\subsubhead\nofrills{b1) } \endsubsubhead
$X$ is contained in a variety of dimension $n+1$,
$V_1$, containing $\infty ^{n - 1} \Bbb P^3$. If we cut $V_1$
with a hyperplane, we find a variety of dimension $n$
containing $\infty ^{n - 1}$ planes, the same situation as in
$a1)$:  such varieties exist if $n\geq 4$ and are classified
only for $n=4$.
\subsubhead\nofrills{b2) } \endsubsubhead $deg \ Y\leq n - 1$.
If $n=3$, then $X$ is a $3$--fold of $\Bbb P^5$ containing
$\infty ^2$ quadric surfaces; by cutting with a hyperplane we
see that either $X$ is a cone over a Veronese surface, or a
rational normal scroll.

If $n=4$, then $X$ is a $4$--fold of $\Bbb P^6$ containing a
family $\Cal F$ of dimension $3$ of cubic or quadric surfaces;
also in this case, by cutting with a general hyperplane, we
reduce to the case $a2)$.

If the surfaces are quadrics, we are
able to say a little bit more: since any quadric contains a
family of dimension $1$ of lines, then either $X$ contains a
family of dimension $4$ of lines, or $X$ contains a family of
lines of dimension $3$ such that through the general one there
are $\infty ^1$ quadrics of $\Cal F$. In the former case,
 by the results quoted in \S 2, if $X$ is
smooth, it is either a scroll in $\Bbb P^2$'s over a surface
or a quadric bundle over a smooth curve: it is easy to see
that both possibilities cannot happen. In the second case, there
are $\infty ^1$ quadrics and a finite number of lines through a
general  point $p$ of $X$: if the quadrics were smooth, for any
$p$ all the quadrics through $p$ should be tangent at $p$,
which is impossible; so the quadrics of the family are cones.

\head 4. $N - n=3$ \endhead
In this section $X \subset \Bbb P^{n+3}$ is a variety
containing a family $\Cal F$ of dimension $h+1$ of
subvarieties $Y$ of dimension $n - h$, spanning a $\Bbb
P^{n-h+2}$. By (1.3) there are the following three
possibilities:
\roster
\item $X \subset V_1$, with $dim \ V_1=n+2$ and $V_1$ contains
a family of dimension $h+1$ of $\Bbb
P^{n-h+2}$;
\item $X \subset V_2$, with $dim \ V_2=n+1$ and $V_2$ contains
a family of dimension $h+1$ of subvarieties $Z$ of
dimension $n-h+1$ and degree $h+1$ each one generating a $\Bbb
P^{n-h+2}$;   \item $deg \ Y$ is bounded by a function of
$h$; by 1.5, (ii) $deg \ Y\leq (h+1)(3h+2)/2$. \endroster

\smallpagebreak
$\boxed{a) \  h=n-1}$

 The varieties $Y$ are
curves of $\Bbb P^3$.
\subsubhead\nofrills{a1) } \endsubsubhead  In case ($1$) the
discussion is analogous to $a1)$ of $\S 3
$; there are no varieties $V_1$ if $n \leq
 2$, they are classified if $n=3$.
\subsubhead\nofrills{a2) } \endsubsubhead In case (2) $V_2$
contains $\infty ^n$ surfaces of $\Bbb P^3$: this is case $b2)$
of $\S 3$.
\subsubhead\nofrills{a3) } \endsubsubhead In case ($3$)
$deg \ Y\leq n(3n-1)/2$.

Let $n=2$: $X$ is a surface of $\Bbb P^5$ containing $\infty
^2$ curves of $\Bbb P^3$; case $a1)$ is impossible, so the first
order foci on a fixed space $H$ of the family $\Phi _1$ form a
quadric; it can be thought of as the union of the focal
lines, i. e. of the intersections of $H$ with the
spaces of $\Phi _1$ which are $\lq\lq$ infinitely near'' to
$H$.  In case $a2)$ $X$ either is contained in a cone over a
Veronese surface or in a rational normal scroll. In case $a3)$,
$deg \ Y\leq 5$; the second order foci
$F_H^2$ on a general space $H$ of the family $\Phi _1$ form a
curve of degree $5$ and genus $2$ which is linked to a line on
the focal quadric.

In the papers \cite{S1} and \cite{S2} C.Segre gave a
classification of the surfaces $X$ as in \it a3)\rm; but there
are some gaps and some rather obscure points in his proof.
Here we give a new proof of his result which follows as much
as possible the original one.
 \proclaim{4.1. Theorem (C.Segre)}
Let $S \subset \Bbb P^5$ be a  surface containing a
family $\Cal F$ of dimension $2$ of irreducible non--degenerate
curves $Y$ of $\Bbb P^3$, parametrized by an irreducible
surface $Z$. If $S$ is of isolated type and not a cone, then one
of the following happens:
\roster
\item $Y$ is a rational cubic
and $S$ is a rational normal scroll of degree $4$;
\item $Y$ is
an elliptic quartic and $S$ is an elliptic normal scroll of
degree $6$;
\item $S$ is a rational surface, isomorphic to a
blowing--up of $\Bbb P^2$, embedded in $\Bbb P^5$ by a linear
system of cubics;
\item $S$ is a rational surface isomorphic to
a blowing--up of $\Bbb P ^2$ embedded in $\Bbb P^5$ by a linear
system of quartics, such that the images of the lines are
rational quartics of $\Bbb P^3$.
\endroster
\endproclaim
\demo{Proof} We give the proof assuming first that $S$ is
smooth. At the end we will indicate how the argument can be
modified if $S$ is singular.

Let us remark first that there are
no base points  for the family $\Cal F$; otherwise, by
projecting from a base point $q$,  we find a surface $S'$ of
$\Bbb P^4$ containing $\infty ^2$ plane curves: $S'$ is a
Veronese surface or a rational normal scroll of degree $3$; so,
since $S$ lies on the cone of vertex $p$ over $S'$, it  is not
of isolated type.

Then, observe that a general curve $Y$ is
smooth; otherwise, if there is a point $p$ which is singular
for all curves of $\Cal F$,  by projecting $S$ from $p$ we
would find as above that $S$ is not of isolated type. If the
singular points of the curves of $\Cal F$ are variable, let $p$
be a variable singular point of $Y$: it should lie on all the
focal lines of the space $H$ which is spanned by $Y$,
so for any $Y$ the focal quadric on $H$ should be a quadric
cone and $p$ its vertex. There are only two possible cases:
either $Y$ is a quartic with a double point or it is a quintic
with a triple point; in both cases $Y$ is rational so the
surface $S$ is rational too. Since the singularities of $Y$
are variable, then $Y$ varies in a linear system of dimension
at least $3$ whose general curve is smooth.

 The intersection number
$j=Y^2$ of two curves of $\Cal F$ is constant; since it may
be computed by intersecting $Y$ with a curve $Y'$ of a \lq\lq
infinitely near'' space of $\Phi _1$, by monodromy we get
that it is equal to the number of second order foci
 of $Y$ on a focal line, i. e. to the number of variable
intersections of $Y$ with the focal lines.  We  have the
following possibilities:
\roster
\item"(a)" $j=1$, $p_a(Y)=0$, deg$Y=4$ or $3$;
\item"(b)" $j=2$, $p_a(Y)=2$, deg$Y=5$;
\item"(c)" $j=2$, $p_a(Y)=1$, deg$Y=4$;
\item"(d)" $j=2$, $p_a(Y)=0$, deg$Y=3$.
\endroster
In case (a), the family $\Cal F$ is linear (a homaloidal
net); it defines a morphism $\phi :S \to \Bbb P^2$ which
sends the curves of $\Cal F$ to the lines, so (3) or (4)
happens.

Let us assume now that $j=2$.\smallpagebreak

\it Claim (4.2).  Let $p, q$ be two
general points of $S$; then there are at least two curves of
$\Cal F$ through them. If $p_a(Y)>0$, there are exactly two.

\it Proof.\rm Assume that the
claim is not true; let $p',q'$ be the intersection of two
general curves of $\Cal F$: then there are infinitely many
curves of $\Cal F$ through $p'$ and $q'$, so they are
fundamental points of $\Cal F$ and by consequence of $\Phi
_1$; the infinitely many linear spaces containing $p'$ and
$q'$ contain  the whole
 line $\langle
p', q' \rangle$, so it consists of fundamental points of
$\Phi _1$ and is focal of first and second order, against the
assumption that $S$ is  of isolated type.

Assume now that there are $r>2$ curves of $\Cal F$ through
two general points of $S$. Let us fix a general point $p$ on
$S$ and denote by $\Cal F _p$ the subfamily of $\Cal F$ of
curves through $p$; it is parametrized by a curve $C_p\subset
Z$. $\Cal F_p$ is a family of curves of dimension $1$, degree
$1$ and index $r$ (let us recall that the \it degree \rm of
$\Cal F_p$ is the number of variable intersections of two
curves of $\Cal F_p$, and the \it index \rm of $\Cal F_p$ is the
number of curves of $\Cal F_p$ passing through a general point
of $S$); in fact, two general curves of $\Cal F _p$ intersect
at a point $q$ different from $p$ and there are $r$ curves of
$\Cal F_p$ through a general $q$. So, if $Y\in\Cal F_p$, we may
define a rational map $g_Y$ from $C_p$ to $Y$ sending a curve
$C$ to $\lbrace Y\cap C\rbrace -p$. The fibers $\lbrace
g_Y^{-1}(q)\rbrace _{q\in Y}$ form a pencil $\Pi _Y$ of
divisors on $C_p$. There is an algebraic family of such
pencils, one for any $Y$ in $\Cal F_p$; note that if $Y\neq Y'$
are curves of $\Cal F_p$, then obviously $\Pi _Y\neq\Pi_{Y'}$.
The classical theorem of Castelnuovo--Humbert asserts that
there cannot be an algebraic family of irrational pencils on a
curve (see \cite{M} for a modern version); it follows that the
curves of $\Cal F$ are rational. This concludes the proof of
Claim (4.2).\smallpagebreak

Let $S$ be as in (d). The curves of $\Cal F$ are skew cubics
with selfintersection $2$. Observe that in this case $S$ is
rational. Let us consider the
exact sequence $$0\longrightarrow \Cal
O_S\longrightarrow\Cal O_S(Y)\longrightarrow\Cal
O_Y(Y)\longrightarrow0;$$ we have $h^0(\Cal O_S(Y))=4$ which
means that $S$ contains in fact a linear family of dimension
$3$ of skew cubics. In particular there are $\infty ^2$ skew
cubics through any point $p$ of $S$; if we project $S$ from
$p$, we get a surface $S'$ in $\Bbb P^4$ with a
$2-$dimensional family of conics, all intersecting a line,
the image of (the tangent plane at) $p$. So it is clear that
both $S'$ and $S$ are rational normal scrolls.

Let $S$ be as in (c); fix $q$, a general point  on $S$, and
$Y$ in $\Cal F_q$: there is a birational map  $Y-
-\to C_q$ (the curve in $Z$ parametrizing $\Cal
F_q$), which associates to $p\in Y $ the point of $C_q$
corresponding to the unique curve of $\Cal F_q$ through $p$. So
$C_q$ is an elliptic curve with the same modulus as all the
curves in $\Cal F$.

We will construct now some rational curves on $S$: let us fix
a $g_2^1$, i.e. a linear series of degree $2$ and dimension
$1$, on $C_q$; it gives a rational family of pairs of curves
of $\Cal F$ through $q$. The variable intersections of such
pairs of curves form a rational curve $L$ on $S$; as the
$g_2^1$ varies on $C_q$, we get a family $\Cal L$ of rational
curves, such that there is exactly one curve of $\Cal L$
through a general point of $S$. Note that $\Cal L$ is
an elliptic pencil, with the same modulus as any $Y$ of $\Cal
F$; in fact we may construct a rational map associating to
$p\in Y$ ($Y\in \Cal F_q$) the curve of $\Cal L$ through $p$,
which does not intersect $Y$ elsewhere. This implies also that
$Y\cdot L=1$ because for $L$ general the intersection $Y\cap
L$ is transversal. Since there is a curve $Y$ of $\Cal F$
through any two points of $L$, such a  $Y$ must split:
$Y=L\cup\Lambda$; so there is a subfamily of dimension $1$ of
$\Cal F$ whose curves are of the form $L\cup\Lambda$, where $L$
is rational, $\Lambda$ is elliptic, $L\cdot \Lambda=1$ and
$\Lambda ^2=0$. Since $deg \ Y=4$, then $\Lambda$ is a plane
cubic and $L$ is a line. Therefore $S$ is an elliptic  scroll,
having elliptic quartics pairwise meeting in two points as
unisecant curves. We conclude that $S$ is an elliptic normal
scroll  of degree $6$.

It remains to exclude the case (b). Note first that in this
case $S$ is birational to the symmetric product $Y^{(2)}$ of
$Y$, where $Y$ is a general curve of $\Cal F$; in fact, as in
case (c), fixed $p\in S$, the system $\Cal F_p$ is
parametrized by $Y$ and we may associate to a point $q$ of $S$
the pair of curves of $\Cal F_p$ passing through $q$. So $S$
is  birational to the Jacobian of a curve of genus
$2$.\smallpagebreak

\it Claim (4.3).  $S$ does not contain
any rational curve; in particular, if it is smooth, it is
minimal.

 Proof. \rm  Note first that if $D$ is any irreducible curve
on $S$, then $D\cdot Y>0$: in fact if $D\cdot Y=0$, then the
curves
 of $\Cal F_q$, $q\in D$,
 split as $D+\Delta$, where $\Delta$ varies in a
$1$--dimensional family; since $S$ does not contain any
algebraic family of rational or elliptic curves, then $\Delta$
is a plane curve of degree $4$. So we have a $1$--dimensional
family of planes; each of them intersects the spaces of $\Phi
_1$ along a line;  these lines generate a rational normal
scroll which contains $S$, against the assumption that $S$ is
of isolated type.

Let now $D\subset
S$ be a rational curve, with $D\cdot Y=m>0$: we
have a map $\phi:C_p \to D^{(m)}$, the $m^{th}$ symmetric
power of $D$, such that $\phi (Y)=Y\cap D$. Since the family
$\Cal F_p$ has index two, then the curve $\phi (C_p)$ is
rational; but, if $m>1$, $\phi$ should be injective, because
$j=2$, so $m=1$. If we take two points of $D$, there is a curve
$Y$ of $\Cal F$ through them, so $Y$ splits: $Y=D+ \Delta$
where $Y\cdot \Delta=1$ and $\Delta$ is irreducible. There are
the maps $\Phi_D: C_p\to D$, $\Phi_{\Delta}: C_p\to\Delta$
defined by the intersection; $\Phi_D$  is $2:1$ so its fibers
form the unique $g^1_2$ on $C_p$; if $\Delta$ is rational,
also $\Phi_{\Delta}$ is $2:1$, therefore $\Phi_D=\Phi_{\Delta}$,
$D=\Delta$ and $Y=2D$: this is impossible because $deg \ Y=5$,
so $p_a(\Delta)>0$. Since $\Delta$ cannot have genus one
because the family $\Cal F$ has constant moduli, then
$p_a(\Delta)=2$ and $\Delta$ is a plane quartic with
$D\cdot\Delta=1$. This implies that $D$ is a line with
$D^2=0$: the claim follows because $S$ is not ruled.

\smallpagebreak

As a first consequence, note that, by degree reasons, all the
curves $Y$ in $\Cal F$ are irreducible. Moreover, any $Y$ is
an ample divisor on $S$, by the Nakai--Moishezon
criterion; in fact $Y^2=2>0$ and, for any
irreducible curve $C$ on $S$, either $C\cdot Y>0$ or $C$ is a
component of $Y$, i.e. $C=Y$. In particular, we may apply the
 Kodaira vanishing theorem, getting
the relations: $h^1(\Cal O_S(Y))=0$, $h^0(\Cal O_S(Y))
=1$.\smallpagebreak

\it Claim (4.4). $S$ is linearly normal.

 Proof. \rm If $S$ were not
linearly normal, then it would be the isomorphic projection of
a surface $F$, $F\subset \Bbb P^6$; $F$ too should contain
$\infty ^2$ quintics of genus $2$, generating spaces that
intersect two by two along a line. Since these spaces cannot
have a common point, they should intersect a fixed $\Bbb P^3$
along planes and the lines of intersection should lie in this
$\Bbb P^3$. By projecting in $\Bbb P^5$ we would find that
the spaces of $\Phi _1$ intersect two by two along lines of a
fixed $\Bbb P^3$; in particular, all the focal quadrics and,
by consequence, the surface $S$ should lie in
this $\Bbb P^3$: a
contradiction.\smallpagebreak

For any $Y$ we consider now the linear system of the curves
$L$ which are residual to $Y$ in a hyperplane section of $S$:
it is a complete linear system of dimension $1$. In fact, by
the cohomology of the exact sequence $$0 \to \Cal I_S(1) \to
\Cal I_Y(1) \to \Cal I_{Y\mid S}(1) \to 0,$$ we get $h^0(\Cal
I_Y(1))=h^0(\Cal     I_{Y\mid S}(1))=2$, because $h^1(\Cal
I_S(1))=0$ by (4.4).

Let $\Cal G$ be the family described by the curves $L$ as
$Y$ varies in $\Cal F$.\smallpagebreak

\it Claim (4.5).  $Dim \ \Cal G=3$.

 Proof. \rm By
the exact sequence  $$0 \to \Cal I_S(1) \to
\Cal I_L(1) \to \Cal I_{L\mid S}(1) \to 0,$$ we get  $h^0(\Cal
I_L(1))=h^0(\Cal     I_{L\mid S}(1))=1$, because $h^0(\Cal
O_S(Y))=1$. Let us consider in $\Cal F\times\Cal G$ the
correspondence $$ \CD  \lbrace (Y,L)| Y+L \in\Bbb P(\Cal O_S(1))
\rbrace   @>p_1>> \Cal F  \\ @VVp_2V   @.   \\ \Cal G   @.
\endCD $$
 where $p_1$ and $p_2$ are the projections. Since the fibers
of $p_1$ have dimension $1$ and the fibers of $p_2$ have
dimension 0, then $dim \ \Cal G=3$. \smallpagebreak

Let us compute the intersection number $Y\cdot L$. Fix $Y' \in
\Cal F$, $H$ a hyperplane containing $Y'$ and let $S\cap
H=Y'\cup L$. Then $Y\cdot (S\cap H)=5=Y\cdot L+Y\cdot Y'$, so
$Y\cdot L=3$.

All the curves of $\Cal G$ are irreducible: if some $L$ is
 reducible, then there is an
irreducible component $\Lambda$ of $L$ such that $\Lambda\cdot
Y<2$, so $\Lambda$ is a component of any curve $Y$ of $\Cal F$
passing through two of its points; i.e. $\Lambda =Y$; but this
cannot happen because $Y^2=2$. Moreover $L$ is ample; in fact
we  will prove  that $L^2>0$ and $L$ is strictly numerically
effective. $S$ being abelian, we have $\Cal O_S\simeq\omega _S,
\Cal O_L(L)\simeq\omega _L$, so, by the exact sequence
$$0\to\Cal O_S\to\Cal O_S(L)\to\Cal O_L(L)\to 0,$$
 we get $1\leq
p_a(L)\leq 3$. Note that $deg \ L=L\cdot (Y+L)=3+L^2$; if
$L^2=0$, then $L$ is a plane cubic, but in this case $S$ should
contain a family of dimension $3$ of plane cubics, which is
impossible. Let
$D$ be an irreducible curve; since all $L$ are irreducible, then
$D\cdot L>0$, because the curves of $\Cal G$ fill up $S$.  By
the Nakai--Moishezon criterion, we conclude that $L$ is ample.

By Kodaira vanishing we get: $p_a(L)=3$, $L^2=4$, $deg \ L=7$,
$deg \ S=12$.

Let us consider now the intersection $Y\cap L$, where $Y\cup
L=S\cap H$ is a hyperplane section of $S$.
\smallpagebreak

\it Claim (4.6).  $Y\cap L$
consists of three points on a line.

 Proof. \rm Let $H$ be
general in the set of the hyperplanes containing the curves of
$\Cal F$. We know that $Y\cdot L=3$; the points of $Y\cap L$
are singular for $S\cap H$, so either they are singular points
of $S$ or $H$ is tangent to $S$ at these points. The former
case is impossible. In the second case, we have a family of
dimension $3$, $\Cal H$, of hyperplanes which are tangent to
$S$; the set of foci of $\Cal H$ on a general $H$ is a line
$l$ (the \lq\lq characteristic line''), because it is the
degeneracy locus of a map $$\Cal O_H^3\simeq T_{H,\Cal
H}\otimes \Cal O_H \longrightarrow \Cal O_H(1)\simeq \Cal
N_{H\mid\Bbb P^5}.$$
 It is well known (\cite{S}) that the tangency points lie on
the characteristic line, so we have the   claim.
\smallpagebreak

Hence, for any hyperplane containing $Y$, we have a
trisecant line of $Y$ (which is a line of the focal
quadric).

As a consequence of (4.6), we have that, for fixed $L$ in
$\Cal G$, the intersection $L'\cap L$ is fixed if $L'\in
\mid L\mid$. In fact, the curves of $\mid L\mid$ are all
residual to the same $Y$; let $\pi _Y=\langle Y\rangle$: since
$L$ is contained in a hyperplane containing $Y$ then the
intersection $\pi _Y\cap L$ consists of $7$ points; moreover
$S\cap\pi _Y=(L\cup Y)\cap \pi _Y= (L\cap\pi _Y)\cup Y$, so
$L\cap\pi _Y$ depends only on $Y$, and is equal to $L\cap
L'$ for any $L' \in\mid L\mid$. None of these points lies on
$Y$; otherwise all the curves of $\mid L\mid$ should have
the same intersection with $Y$ (for any point on $Y$ there is
only one trisecant), so a curve in $\mid L\mid$ passing
through another point of $Y$ should have $Y$ as a component,
which has been excluded.

To conclude the proof, we will finally consider the family of
the trisecants of $S$. Let $t$ be such a line and $Y$ be a
curve of $\Cal F$ passing through two of the intersection
points. Then $t \subset\pi _Y$, so either the $3^{rd}$
intersection point lies on $Y$ and $t$ is a trisecant of $Y$, or
it is one of the base points of $\mid L\mid$. Trisecants of the
first type form a family $\Cal T$ of trisecants of $S$ of
dimension $3$; those of the second type a family of dimension
$2$: from now on, we will consider only the trisecants of
$\Cal T$.

Fix $p\in S$: the trisecants through $p$ are parametrized by
$C_p$, because there is exactly one trisecant for any
curve of $\Cal F$ through $p$, so they form a cone $K_p$ over
a curve of genus $2$, whose degree is at least $4$. If $deg \
K_p=4$, then $K_p$ is a cone over a plane quartic, so $\langle
K_p\rangle$ has dimension $3$; the intersection points of the
trisecants of $K_p$ generate a curve contained in $S$ into this
$3$--space, whose degree is at least $8$. Since curves
corresponding to different points of $S$ are different, $S$
should contain also $\infty ^2$ curves of $\Bbb P^3$ of degree
at least $8$: a contradiction. Hence, $deg \ K_p\geq 5$; in
particular the number of trisecants through $p$ contained in a
hyperplane, counted with multiplicities, is at least $5$.

Let us
fix now a curve $Y$ of $\Cal F$ not containing $p$ and consider
$H=\langle Y,p\rangle$, $H\cap S=Y\cup L$, $L\in \Cal G$, $p\in
L$. We will see now that there are at most $4$ trisecants of
$S$
contained in $H$ and passing through $p$: this gives a
contradiction that concludes the proof of the theorem.

Let $t$ be such a trisecant; since $p\notin \langle Y\rangle$,
then deg$(t\cap Y)\leq 1$, otherwise $t\subset \langle
Y\rangle$. So either $t$ is a trisecant of $L$ or it is a
secant of $L$ intersecting $Y$. We prove that if $p$ is  a
general point of $L$, there is at most one trisecant of $L$
through $p$. If there are two trisecants through $p$, $t_1$ and
$t_2$, let us consider the projection from $t_i$,
$\pi _i$, $i=1,2$. If it is birational, then $\pi _i(L)$  is a
plane curve of degree $4$ with a node, hence of genus $2$;
which is impossible. So $deg \ \pi _i=2$. Note that $L$ is
hyperelliptic: the $
g^1_2$ is cut out by the $\Bbb P^3$'s through the plane
$\langle
t_1,t_2\rangle$; the pairs of points having the same image via
$\pi _i$ form the unique $g^1_2$: since they are contained
both in planes through $t_1$ and through $t_2$, then they are
contained in lines passing through $p$. It is clear that this
cannot happen for a general point of $L$.

 On the other hand, to compute the number (with multiplicity)
of secants of $L$ through $p$ intersecting $Y$, let us assume
that $p$ is one of the $4$ base points of $\mid L\mid$
belonging to $\pi _Y-Y$. In this case, such a trisecant $t$
lies inside $\pi _Y$ and there are the following possibilities:

(i) $t$ is one of the $3$ lines joining $p$ with a point of
$L\cap Y$;

(ii) $t$ passes through $2$ base points of $\mid L\mid$ and
meets $Y$.

In case (i), $t$ has to be counted with multiplicity one, if
$L$ is general, otherwise it should  meet $Y$ elsewhere; but in
the linear system $\mid L\mid $ there is only  a finite number
of curves passing through the intersection points of the
secants of $Y$ through $p$. Possibility (ii) does not happen,
if $Y$ is general in $\Cal F$: in fact, assume that, for any
$Y$ of $\Cal F$, a line through $2$ of the base points of $\mid
L\mid$ intersects $Y$ and let $x_Y$ be the intersection point.
Let $\bar L$ be the curve of $\mid L\mid$ through $x_Y$: the
other points $x_1, x_2$, of $Y\cup\bar L$ are on a line through
$x_Y$, so there are $2$ trisecants through $x_Y$. Arguing as
above (see \cite{S2}), we get that $x_Y$ lies on all lines
generated by the pairs of points of the $g^1_2$ of $L$; so
$x_Y$ is the vertex of a cone of trisecants $S$ having degree
$3$. Since, as $Y$ varies in $\Cal F$, the points $x_Y$ cover
$S$, we have the required contradiction.

If $S$ is singular, we may consider $\tilde S$, a minimal
desingularization of $S$ and repeat the above argument on
$\tilde S$. There are two points where one must be careful.
Precisely, concerning the proof of Claim (4.3), if $D$ is an
irreducible curve on $\tilde S$ such that $D\cdot \tilde Y=0$
(where $\tilde Y,$ varying in $\tilde \Cal F$, denotes the
preimage on $\tilde S$ of a curve $Y$ of $\Cal F$), then
there is a subfamily of $\tilde \Cal F$ of curves of the
form $D+\Delta$, where $\Delta$ varies in a $1$-dimensional
family. In addition to the above discussed case, it may
happen that $D$ is a line which goes to a singular point of
$S$ and $\Delta$ is a curve of genus $2$ whose image on $S$
is a quintic of the family $\Cal F$. But in this case $D\cdot
\Delta=D\cdot \tilde Y=1$ against the assumption. So $\tilde
S$ may contain some rational curves, precisely lines $D$ with
$D\cdot\tilde Y=1$, but they contract to isolated singular
points on $S$. Moreover, in the proof of Claim (4.6), we
have to exclude that, for $H$ general, some of the $3$
points of $Y\cap L$ lies on the double curve $D'$ of $S$.
 If $D'$ is a
double line, keeping a hyperplane containing $Y$ and
$D'$, we would find a non-integral residual to $Y$, which is
impossible. So $D'$ should have degree $2$ or $3$ and its
points of intersection with $Y$ should be base points for
the linear system $\mid L\mid$ of the residuals. If there are
$3$ points, by imposing  the passage through a $4^{th}$
point of $Y$, we would find a reducible $L$; if there are
$2$, then the third intersection of the curves of $\mid
L\mid$ with $Y$ would give a rational parametrization of
$Y$: therefore both cases are excluded.

 \enddemo

\remark{4.7.
Remarks}

 (i) Let $S$ be a
surface of type (4) in (4.1); to get such a surface one has to
find a linear system $\delta$ of dimension $5$ of plane
quartics such that, for any line $l$, the curves of $\delta$
containing $l$ form a pencil. If one wants a linearly normal
surface, then $\delta$ has the form $\delta =\mid
4l-\sum_{i=1}^r P_i\mid$, where some of the points $P_i$ may
coincide or be infinitely near; so the cubic component of the
quartics containing a line $l$ varies in the linear system
$\delta (-l)=\mid 3l-\sum_{i=1}^r P_i\mid$. We sketch here
how one can proceed provided the points $P_i$ are pairwise
distinct. Note that $r\leq 9$ distinct points impose dependent
conditions to the quartics if and only if $6$ of them lie on a
line. We exclude this case, because then the line is a fixed
component of $\delta$, and we fall in case (c).

Let $r=9$; $S$ is a surface of degree $7$; dim$\delta (-l)=1$
if and only if $P_1,...,P_9$ impose less than $9$ conditions
to the cubics, i.e. either they are the complete intersection
of two cubics, or $8$ of them lie on a conic. In the first
case, as $l$ varies $\delta (-l)$ is a fixed pencil of cubics,
so $\delta =\Bbb P(V)$, where $V\subset H^0(\Cal O_{\Bbb
P^2}(4))$ has the form $V=(x_0,x_1,x_2)(F,G)$ ($x_0,x_1,x_2$
homogeneous coordinates in the plane, $F,G$ polynomials of
degree $3$ defining two cubics of the pencil). This means that
$S\subset\Bbb P^1\times\Bbb P^2$, i.e. it is not of isolated
type. In the second case, $\delta (-l)$ is also a fixed pencil
of cubics, having a common conic $\Gamma$; $\delta$  does not
separate the points of $\Gamma$, so the surface $S$ is
singular: it has a quadruple point $p$ image of $\Gamma$. If we
project $S$ in $\Bbb P^4$ from $p$, we get a surface $S'$ of
degree $3$, i.e. a rational normal scroll. Since $S$ is
contained in the cone of vertex $p$ over $S'$ it is not of
isolated type.

 Coming to the non--linearly normal surfaces, $S$ is
projection of a surface $F$ of $\Bbb P^n$, $n\geq 6$, which is
also embedded by a linear system of quartics. There are two
cases: either the rational quartics of $F$ are already in $\Bbb
P^3$ and the projection is general, or they
are rational normal quartics and one has to find a suitable
projection. In the first case, $F$ is embedded in $\Bbb P^n$ by
a linear system $\delta =\mid 4l-\sum _{i=1}^{14-n} P_i\mid$
such that $\delta (-l)=\mid 3l-\sum _{i=1}^{14-n} P_i\mid$ has
dimension $n-4$. As before we exclude that the points ${P_i}$
lie on a line. So the only possibility is $n=6$; then
${P_1,...,P_8}$ are on a conic $\Gamma$ and $\delta$ does not
separate the points of $\Gamma$. $F$ has a quadruple point $P$
image of $\Gamma$; projecting from $P$, we find the Veronese
surface of $\Bbb P^5$. Also in this case $S$ is not of
isolated type, because it projects on the Veronese surface in
$\Bbb P^4$.

 The second case is more complicated:
$\delta =\mid 4l-\sum _{i=1}^{14-n} P_i\mid$ is a projective
space of dimension $n$; for any line $l$ in the plane, $\delta
(-l)$ is a linear subspace of dimension $n-5$ of $\delta$. We
have to find a subspace $\delta '\subset \delta$, of dimension
$5$, such that dim$(\delta '\cap\delta (-l))=1$ for any line
$l$.

If $n=6$, there are $8$ base points, the only possible choice
for $\delta '$ is imposing the passage through the $9^{th}$
base point of the pencil of cubics through $P_1,...,P_8$: we
fall in  one of the previous cases.

Assume $n=7$: for any line $l$, $\delta (-l)$ is a plane
into $\delta$. Let us consider $X=\cup \delta (-l)$: it is a
variety of dimension $4$ in $\Bbb P^7$, ruled by planes. We
may identify the net of cubics through $P_1,...,P_7$ with
$\Bbb P^2$; $\check \Bbb P^2\times\Bbb P^2$, the set of
pairs $(l,\Gamma$), $l$ line of $\Bbb P^3$, $\Gamma$ cubic
through $P_1,...,P_7$, is embedded in $\Bbb P^8$ via the
Segre map, as a variety of degree $6$, and $X$ is a
projection of its, so $deg \ X\leq 6$. We would like to
find a subvariety $Y\subset X$, ruled by lines, of dimension
$3$, generating a $\Bbb P^5$, so $Y$ should have degree at most
$5$. For any hyperplane $H$ containing $Y$, we would find a
variety $Y'$ of dimension $3$, residual to $Y$ in $H\cap X$, and
such varieties $Y'$ should describe a linear system of dimension
$1$ inside X. If $deg \ Y=5$, then the $Y'$ should be linear
spaces, which is impossible. If $deg \ Y=3$, $Y$ is the rational
normal scroll $\Bbb P^1\times\Bbb P^2$; this means that the
pencil of cubics is fixed, i.e. the projection is centered in
two points of $F$ and we find again the above cases. If
$deg \ Y=4$, then its general hyperplane section is a surface of
degree $4$ in $\Bbb P^4$, therefore a Veronese surface, or a
Del Pezzo surface, complete intersection of two quadrics,
or a cone, or a
rational non--normal scroll. In none of the former three cases
$Y$ can be a scroll in lines; in the last case, $Y$ is a
rational non-normal $3-$fold ruled by planes: these planes
should belong to one of the rulings of $X$, which is
impossible.

It is clear that, if $n>7$, a similar analysis becomes more and
more complicated and we were not able to conclude
it.\smallpagebreak

(ii) It easy to see that the rational normal scroll of degree
$4$ (case (1) of the Theorem) is the unique case of a surface
of $\Bbb P^5$ with a $3$--dimensional family of curves of $\Bbb
P^3$. In fact, given such a surface $S$ and a point $p$ on $S$,
the curves of the family through $p$ form a subfamily of
dimension $2$; if we project $S$ from $p$, we get a surface $S'$
in $\Bbb P^4$ with a $2$--dimensional family of plane curves,
all intersecting a line and of degree $2$. So it is clear that
$S'$ is a rational normal scroll.\smallpagebreak

(iii) If there exists a surface $S$  of isolated
type containing $\infty ^2$ rational quartics, then
for any point $p$ of $S$ and any curve $Y$ of $\Cal F$ through
$p$, there is a finite number of trisecants to $Y$;
so there is a cone of trisecants $S$ at any point $p$. Hence
such a surface would provide an example of a surface of $\Bbb
P^5$ with  a family of dimension $3$ of  trisecant lines,
filling up a variety of dimension $4$. Up to now, the only
other known example is a non-general Enriques surface of
degree $10$ (see \cite{CV}).
\endremark
\smallpagebreak
$\boxed{b) \  h=n-2}$

The varieties $Y$ are
surfaces of $\Bbb P^4$.  If $n=3$, then $h=1$ and
$X$ is a threefold of $\Bbb P^6$ containing a family of
dimension $2$ of surfaces of $\Bbb P^4$. There are as usual $3$
possibilities: \subsubhead\nofrills{b1) } \endsubsubhead $X$ is
contained in a variety of dimension $5$ containing $\infty ^2$
spaces of dimension $4$: but this cannot happen for $X$
non--degenerate
 by  $\S 2$;
\subsubhead\nofrills{b2) } \endsubsubhead $X$ is contained in
a cone of dimension $4$ over a Veronese surface or in a
rational normal scroll of dimension $4$;
\subsubhead\nofrills{b3) } \endsubsubhead The second order
foci are Castelnuovo surfaces. In this case a general
hyperplane section of $X$, $X\cap H=S$, is a surface of $\Bbb
P^5$ of isolated type containing a family of dimension $2$ of
surfaces of $\Bbb P^3$, hence a surface  as in Theorem 4.1. If
$S$ is a rational normal scroll, then $X$ exists and is a
rational normal scroll of degree $4$ too.

If $S$ is as in (4),
then $X$ has to contain a family of surfaces of degree $4$
with rational sections, i.e. of Veronese surfaces, or cones,
or rational non-normal scrolls (with a node). Each of them
must lie in the variety of $1^{st}$ order foci, which is a
quadric, so the  case of the Veronese surfaces cannot happen.
If the surfaces are all singular
 with a fixed singular point $p$, either $X$ is a cone or the
projection from $p$ gives a $3$-fold $X'$ in $\Bbb P^5$
containing $\infty ^2$ surfaces of $\Bbb P^3$ of degree at
most $2$; so $X'$ is not of isolated type: it should be
contained in a $4$-fold of $\Bbb P^5$ containing a
$2$-dimensional family of $\Bbb P^3$'s, but this is impossible
(case $b1)$). If $p$ is a variable singular point of the
surfaces, arguing as in 4.1, we see that it is the vertex of
the focal quadric cone and that this vertex is fixed: a
contradiction.

For discussing the other two
cases, let us restrict to smooth threefolds. If $S$ is  an
elliptic normal scroll as in (2), then $X$ has to be linearly
normal in $\Bbb P^6$ of degree $6$: by Ionescu classification
(see \cite{I}), it does not exist. Finally, if $S$ is as in
(3), it is a Del Pezzo surface of degree $5$ linearly normal in
$\Bbb P^5$, or a projection of a Del Pezzo surface of $\Bbb
P^n$, $5<n\leq 9$.  $X$ should have sectional genus $1$, hence
it should be either a threefold ruled by planes over an elliptic
curve (which cannot happen in our case) or a rational
threefold. Such rational threefolds exist (see \cite{Sc}) of
degree $n$, for $n=5,...,8$; they may be realized as images of
$\Bbb P^3$ via some maps given by linear systems of cubic
surfaces having a base curve of degree $9-n$.

 \enddocument